\documentclass[a4paper,twoside,reqno]{amsart}
\usepackage[utf8]{inputenc}
\usepackage{amsthm,amssymb,amsmath,amsfonts,mathrsfs,nccmath}
\usepackage{orcidlink}
\usepackage{mathtools,stmaryrd}
\usepackage{hyperref}
\hypersetup{
  colorlinks=true,
  urlcolor=green,
  linkcolor=blue,
  citecolor=green
}
\newtheorem{thm}{Theorem}[section]

\theoremstyle{definition}
\newtheorem{defn}[thm]{Definition}

\numberwithin{equation}{section}

\newcommand{\ba}{\begin{aligned}}
\newcommand{\ea}{\end{aligned}}
\newcommand{\norm}[1]{\big\Vert#1\big\Vert}
\newcommand{\jump}[1]{\llbracket#1\rrbracket}
\newcommand{\R}{\mathbb R}
\newcommand{\T}{\mathbb{T}}
\newcommand{\p}{\partial}
\newcommand{\rd}{\mathrm{d}}
\newcommand{\dvg}{\nabla\cdot}
\newcommand{\Lap}{\Delta}
\newcommand{\wt}[1]{\widetilde{#1}}

\newcommand{\ga}{\gamma}
\newcommand{\de}{\delta}

\newcommand{\om}{\omega}
\newcommand{\Ga}{\Gamma}
\newcommand{\Om}{\Omega}
\newcommand{\ka}{\kappa}
\newcommand{\veps}{\varepsilon}

\newcommand{\sig}{\sigma}
\newcommand{\bx}{\overline{x}}

\newcommand{\be}{\begin{equation}}
\newcommand{\ee}{\end{equation}}
\newcommand{\bes}{\begin{equation*}}
\newcommand{\ees}{\end{equation*}}

\newcommand{\cT}{\mathcal{T}}
\newcommand{\cU}{\mathcal{U}}
\newcommand{\cR}{\mathcal{R}}

\newcommand{\Ompm}{\Omega^{\pm}}
\newcommand{\ppm}{^{\pm}}
\newcommand{\TT}{\mathbb{T}^2}

\title[From entropy waves to vortex sheets]{From entropy waves to vortex sheets for the Euler equations}

\author[Trakhinin]{Yuri Trakhinin}
\address{Sobolev Institute of Mathematics, Siberian Branch of the Russian Academy of Sciences,
Koptyug av. 4, 630090 Novosibirsk, Russia}
\email[Y. Trakhinin]{trakhin@math.nsc.ru}

\author[WANG]{Wei Wang}
\address{School of Mathematical Sciences,
Zhejiang University, Hangzhou, People’s Republic of China}
\email[W. Wang]{wangw07@zju.edu.cn}

\author[ZHANG]{Zhifei Zhang}
\address{School of Mathematical Sciences, Peking University, Beijing, People’s Republic of China}
\email[Z. Zhang]{zfzhang@math.pku.edu.cn}

\author[ZHAO]{Wenbin Zhao$^{\dagger}$\orcidlink{0000-0003-0804-1934}}
\address{School of Mathematics, Renmin University of China, People’s Republic of China}
\email[W. Zhao]{wenbizhao2@ruc.edu.cn}
\thanks{$\dagger$Corresponding author.}

\date{\today}
\begin{document}
\begin{abstract}
We construct piecewise smooth solutions of the full compressible Euler equations that evolve from entropy waves into vortex sheets while retaining a strictly positive Taylor jump. For every integer $\ka\geq1$, we obtain solutions satisfying initial material compatibility for pressure and velocity through order $\ka$, whose tangential velocity jump has a nonzero leading term of order $t^{\ka+1}$. These solutions are periodic in the tangential variables and admit a smooth realization in a slab with impermeable walls. We also construct solutions on a local nonperiodic tangential patch that satisfy initial material compatibility at every order, yet have a tangential velocity jump that is nonzero for every sufficiently small positive time and flat at $t=0$. The finite-order construction uses a normal Cauchy--Kowalevski problem for velocity potentials, whereas the all-order construction uses a homogeneous reduction to a smooth hyperbolic Cauchy problem.
\end{abstract}
\maketitle

\section{Introduction}
\subsection{Presentation of the problem}
In this article, we consider the full compressible Euler equations
\begin{equation}\label{eq:euler-conservative}
\begin{cases}
\p_t\rho+\dvg(\rho u)=0,\\
\p_t(\rho u)+\dvg(\rho u\otimes u)+\nabla p=0,\\
\p_tE+\dvg\big((E+p)u\big)=0,
\end{cases}
\end{equation}
where $\rho$, $u$, $p$, and $E$ are the density, velocity, pressure, and total energy density, respectively. The equation of state and the energy are
\begin{equation*}
p=A\rho^\ga e^S,\qquad
E=\frac{p}{\ga-1}+\frac12\rho|u|^2,
\end{equation*}
where $S$ is the entropy and $A>0$, $\ga>1$ are constants.

Following the graph convention of \cite{WZZ}, we consider piecewise smooth solutions in the domain $\Om=\TT\times(-L,L)$, $L>0$, with
\begin{equation*}
\Ompm(t)=\{(\bx,x_3)\in\Om:\ \pm(x_3-f(t,\bx))>0\},
\qquad
\Ga_f=\{x_3=f(t,\bx)\}.
\end{equation*}
Here $\T=\R/(2\pi\mathbb Z)$. The upward normals are $N=(-\p_1f,-\p_2f,1)$ and $n=N/|N|$, and we write $\jump h=h^+|_{\Ga_f}-h^-|_{\Ga_f}$.

We study contact discontinuities, for which the relative normal mass flux $\rho\ppm(u\ppm\cdot N-f_t)$ vanishes. The Rankine--Hugoniot conditions then reduce to
\begin{equation}\label{eq:contact}
\jump p=0,\qquad u^+\cdot N=u^-\cdot N=f_t.
\end{equation}
A contact discontinuity is called an entropy wave if $\jump u=0$ and $\jump\rho\neq0$, and consequently $\jump S\neq0$. A nonzero tangential velocity jump gives a vortex sheet on the interface.

Within each phase, for positive density and pressure, \eqref{eq:euler-conservative} takes the form
\begin{equation}\label{eq:euler}
D_t\rho=-\rho\dvg u,\qquad
D_t u=-\rho^{-1}\nabla p,\qquad
D_t p=-\ga p\dvg u.
\end{equation}
Here $D_t\ppm=\p_t+u\ppm\cdot\nabla$ denotes the phasewise material derivative; the phase index is suppressed in \eqref{eq:euler} and whenever it is clear from the context. The entropy $S=\log(p/(A\rho^\ga))$ satisfies $D_tS=0$.

At the fixed walls, we impose $u_3\ppm=0$ on $x_3=\pm L$. 
The initial data are $(\rho_0\ppm,u_0\ppm,p_0\ppm,f_0)$, with continuous pressure and full velocity traces and a nonzero density jump on $\Ga_{f_0}$.

We use the unit-normal Taylor jump
\begin{equation}\label{eq:taylor}
\cT=n\cdot\jump{\nabla p}
=\frac{\jump{\nabla_Np}}{|N|},\qquad \nabla_Np=N\cdot\nabla p.
\end{equation}
Its positivity is equivalent to the Taylor sign condition in \cite{WZZ}. Throughout, $q=p\ppm|_{\Ga_f}>0$ denotes the common pressure trace.

For the constructions we write $(x,y)=(\bx,x_3)$ and $u=(U,V)$, and allow $x\in\T^{d-1}$, $d\geq2$. The smooth construction instead uses a local nonperiodic coordinate $x>0$; it embeds in three dimensions by adjoining a passive tangential variable. In two dimensions, the unit tangent oriented toward increasing $x$ is $\tau=(1,f_x)/\sqrt{1+f_x^2}$, and $u_\tau=u\cdot\tau$.

\subsection{Related work and the question of vortex-sheet generation}
The stability of characteristic discontinuities is a central problem in multidimensional conservation laws. Entropy waves and vortex sheets give rise to free boundary problems with different transmission conditions; see the survey of Chen and Wang \cite{MR3289359}. 
For a normal-mode analysis of planar nonisentropic contacts and the distinction between neutral stability and violent instability in two and three space dimensions, we refer to Coulombel and Morando \cite{MR2159807} (see also references therein for the original works).

For two-dimensional isentropic vortex sheets, Coulombel and Secchi established linear estimates \cite{MR2095445} and local nonlinear existence \cite{MR2423311} under a supersonic condition. Morando and Trebeschi \cite{MR2441089} obtained the corresponding linear estimates for nonisentropic flows, and Morando, Trebeschi and Wang \cite{MR3912754} proved nonlinear stability under suitable supersonic assumptions. The nonlinear arguments use Nash--Moser iteration to overcome the loss of derivatives in the linear estimates. These results concern perturbations of existing vortex sheets with a nonzero tangential velocity jump.

Contact waves have also been studied through dissipation. In one space dimension, Huang, Matsumura and Xin \cite{MR2208289} proved asymptotic stability of weak viscous contact waves for the compressible Navier--Stokes equations under a zero excess mass condition. Huang, Xin and Yang \cite{MR2450610} subsequently treated general small perturbations for the Navier--Stokes and Boltzmann equations. For planar entropy waves in three dimensions, Kang, Vasseur and Wang \cite{MR4268832} established stability and uniqueness in a class of zero-dissipation limits of the Navier--Stokes--Fourier system.

In our previous work \cite{WZZ}, we studied entropy waves in Eulerian coordinates. By deriving an evolution equation for tangential derivatives of the interface, we related the Taylor sign condition to its hyperbolicity and obtained a priori estimates without loss of regularity. These estimates concern solutions whose pressure and full velocity remain continuous across the evolving interface. The compatibility conditions needed to preserve this continuity, and the possibility of transition to a vortex sheet when they fail, were also discussed there. The present paper constructs solutions that develop tangential slip from initially continuous velocity and examines how this transition depends on the order of initial material compatibility.

The pressure mechanism behind this question was identified in \cite{WZZ}. If full velocity continuity persists, differentiating along the common material interface trajectories and projecting the momentum equations tangentially gives
\begin{equation}\label{eq:instantaneous-sheet-forcing}
\jump{\rho^{-1}}\nabla_{\Ga_f}q=0,
\end{equation}
where $\nabla_{\Ga_f}$ denotes the surface gradient. A nonzero density jump therefore forces the common pressure trace to be spatially constant on each connected interface component. A tangentially varying trace gives the two phases different tangential pressure accelerations and can generate a vortex sheet.

We consider the hierarchy obtained by differentiating the pressure and full velocity transmission conditions along phasewise material trajectories. We ask whether matching this hierarchy at the initial time, together with a strictly positive Taylor jump, suffices to preserve an entropy wave. We construct contacts satisfying initial material compatibility to any prescribed finite order and, in a separate local construction, to every order. These constructions show that initial material compatibility, even to every order in the smooth setting, together with a positive Taylor jump, does not force a contact discontinuity to remain an entropy wave.

The distinction between an initial entropy-wave configuration and a persisting entropy-wave solution is essential here. Each construction has its own reference solution, obtained by setting $\varepsilon=0$, which remains an entropy wave throughout the local time interval. For $\varepsilon>0$, the constructed solutions have the same interface as the corresponding reference wave and initially satisfy the entropy-wave conditions, yet develop nonzero tangential slip at every sufficiently small positive time while retaining a strictly positive Taylor jump. This occurs despite initial material compatibility through any prescribed finite order in the first construction and at every order in the second.

\section{Main results and ideas of the proof}\label{sec:model}
\subsection{Initial material compatibility}
We first specify the compatibility conditions used in the two constructions.

\begin{defn}\label{def:compat}
Initial phase data have \emph{entropy wave material compatibility of order $\ka$} if, on the initial interface,
\begin{equation}\label{eq:compat}
\jump{(D_t^jp)_0}=0,\qquad \jump{(D_t^ju)_0}=0,
\qquad 0\leq j\leq \ka.
\end{equation}
Every jet is computed recursively from the corresponding phasewise Euler equations. These conditions are necessary for an entropy wave continuation whenever the indicated material derivatives exist classically.
\end{defn}

For a smooth phase field $h\ppm$, write $\underline h\ppm(t,x)=h\ppm(t,x,f(t,x))$, and define $\underline U\ppm$ similarly. The kinematic condition in \eqref{eq:contact} gives
\begin{equation}\label{eq:trace-material}
\left.(D_t\ppm h\ppm)\right|_{\Ga_f}
=\big(\p_t+\underline U\ppm\cdot\nabla_x\big)\underline h\ppm.
\end{equation}
Iteration identifies the material jets with those of the phasewise interface transport operators.

\subsection{Phase laws and the reference entropy wave}
Fix constants $b,P_*>0$ and $r_->r_+>0$. Set $K_\pm=P_*/r_\pm^\ga$ and define
\begin{equation}\label{eq:R-H}
R_\pm(p)=r_\pm(p/P_*)^{1/\ga},\qquad
H_\pm(p)=\frac{\ga P_*}{(\ga-1)r_\pm}
\left((p/P_*)^{(\ga-1)/\ga}-1\right).
\end{equation}
Thus $H_\pm(P_*)=0$, $H_\pm'=R_\pm^{-1}$, and $\rho\ppm=R_\pm(p\ppm)$ is equivalent to $p\ppm=K_\pm(\rho\ppm)^\ga$.

Near $\xi=0$, put
\begin{equation}\label{eq:background}
\bar p_\pm(\xi)=P_*(1-\beta_\pm\xi)^{\ga/(\ga-1)},\qquad
\bar\rho_\pm(\xi)=r_\pm(1-\beta_\pm\xi)^{1/(\ga-1)},
\end{equation}
where $\beta_\pm=(\ga-1)r_\pm b/(\ga P_*)$. These profiles are positive on a small common collar and satisfy
$H_\pm(\bar p_\pm)=-b\xi$ and $\p_\xi\bar p_\pm=-b\bar\rho_\pm$.

With $\xi=y-bt^2/2$, the fields $\rho\ppm=\bar\rho_\pm(\xi)$, $p\ppm=\bar p_\pm(\xi)$, and $u\ppm=bt\,e_y$ form an entropy wave, where $e_y$ is the unit vector in the $y$ direction. Indeed, $D_t\xi=\dvg u=0$, while the pressure gradient supplies the acceleration $b e_y$. This is the reference solution for the finite-order construction.
The reference fields define an exact entropy-wave solution throughout the time interval under consideration: the pressure and the full velocity remain continuous across the moving interface, while the density jump remains nonzero. Thus, the reference solution is a persisting entropy wave, not merely an initial entropy-wave configuration. It is recovered by setting $\varepsilon=0$ in the finite-order construction below.

\subsection{Main results}
We construct two families of contact discontinuities whose initial interfaces are entropy waves. In the first family, the velocity traces separate at a prescribed finite order in time; in the second, their difference is nonzero but flat at the initial time. Both families retain a strictly positive Taylor jump.

For the first construction let
\begin{equation}\label{eq:params}
\ka\in\mathbb N,\quad \ka\geq1,\qquad
a\in C^\om(\T^{d-1};\R),\quad \nabla_xa\not\equiv0.
\end{equation}

\begin{thm}[Finite-order compatibility]\label{thm:main}
Under \eqref{eq:params}, there are $\veps_0,T,\de>0$ such that for every $|\veps|\leq\veps_0$ there is a contact solution of \eqref{eq:euler-conservative} on
\begin{equation*}
0\leq t\leq T,\qquad 0<\pm\big(y-f(t)\big)<\de,
\qquad x\in\T^{d-1},\qquad f(t)=\tfrac12bt^2.
\end{equation*}
The phase fields are smooth up to the interface and depend smoothly on $\veps$. They have the following properties.
\begin{enumerate}
\item The pressure trace is $q(t,x)=P_*+\veps t^\ka a(x)$, the normal velocity is $bt$, and $\rho\ppm=R_\pm(p\ppm)$. The entropy is constant in each phase and differs across the interface.
\item The initial data satisfy \eqref{eq:compat} through order $\ka$. Writing $u=(U,V)$, for $\veps\neq0$ the first nonzero tangential velocity defect occurs at order $\ka+1$ and is given by
\begin{equation}\label{eq:first-defect}
\jump{(D_t^{\ka+1}U)_0}
=-\veps \ka!\left(\frac1{r_+}-\frac1{r_-}\right)\nabla_x a.
\end{equation}
\item On the interface,
\begin{equation}\label{eq:sign-exact}
\cT=b\big(R_-(q)-R_+(q)\big)\geq\theta>0,
\end{equation}
where $\theta$ is independent of $\veps$.
\item For every fixed nonnegative integer $s$, the tangential slip has the expansion
\begin{equation}\label{eq:slip}
\jump{U}(t)
=-\frac{\veps t^{\ka+1}}{\ka+1}
\left(\frac1{r_+}-\frac1{r_-}\right)\nabla_x a
+\cR_\veps(t),\qquad
\norm{\cR_\veps(t)}_{H^s}\leq C_s|\veps|t^{\ka+2}.
\end{equation}
After decreasing $T$ for a fixed $s$, there is $c_s>0$ such that
\begin{equation}\label{eq:lower-bound}
\norm{\jump{U}(t)}_{H^s(\T^{d-1})}
\geq c_s|\veps|t^{\ka+1},\qquad 0<t\leq T.
\end{equation}
\item For every fixed integer $m\geq0$, the initial fields $\cU=(\rho,u,p)$ satisfy
\begin{equation}\label{eq:smallness}
\sum_\pm\norm{\cU_{0,\veps}\ppm-(\bar\rho_\pm,0,\bar p_\pm)}_{H^m(\Om_{0,\de}\ppm)}
\leq C_m|\veps|,
\end{equation}
where $\Om_{0,\de}\ppm=\T^{d-1}\times\{0<\pm y<\de\}$.
\end{enumerate}
For any prescribed slab width $L>0$, the collar can be chosen smaller than $L$ and the initial fields can be smoothly extended to a fixed-wall slab $\T^{d-1}\times(-L,L)$. The extensions yield piecewise $C^\infty$ contacts on a common short time interval with $V=0$ at $y=\pm L$. Estimates \eqref{eq:slip}--\eqref{eq:lower-bound} persist, and \eqref{eq:smallness} holds on the full half-slabs relative to the corresponding smoothly extended reference data. The constants and lifespan may depend on $\ka$.
\end{thm}

The contact solution corresponding to $\varepsilon=0$ remains an entropy wave for $0\le t\le T$. For every sufficiently small $\varepsilon\ne0$, the constructed initial data are arbitrarily close to the reference data in any prescribed finite Sobolev norm and satisfy entropy-wave material compatibility through order $\kappa$. Nevertheless, the tangential velocity jump is not identically zero for every $0<t\le T$. Hence, arbitrarily small perturbations within this class can destroy the entropy-wave property immediately, despite the strictly positive Taylor jump.

\begin{thm}[Smooth infinite-order compatibility]\label{thm:flat-nonlinear}
Fix $\ga>1$, $b,P_*>0$, $r_->r_+>0$, $\veps>0$, and a compact tangential interval $I\subset(0,\infty)$. There are $T,\de>0$ and a piecewise $C^\infty$, phasewise isentropic and irrotational solution of \eqref{eq:euler-conservative} on an open neighborhood of the two collars
\begin{equation*}
0\leq t\leq T,\qquad 0<\pm(y-f(t,x))<\de,
\qquad x\in I,\qquad f(t,x)=\frac{bt^2}{2x},
\end{equation*}
whose interface is a contact discontinuity. Its initial interface traces satisfy
$u_0\ppm=0$, $p_0\ppm=P_*$, $\rho_0\ppm=r_\pm$, and entropy wave material compatibility holds at every order. The Taylor jump is uniformly positive for $0\leq t\leq T$, with
\begin{equation*}
\cT(0,x)=\frac{b(r_--r_+)}{x}.
\end{equation*}
With the unit tangent oriented toward increasing $x$, its slip obeys
\begin{equation*}
\jump{u_\tau}(t,x,f(t,x))
=\veps\left(\frac1{r_+}-\frac1{r_-}\right)
\frac{t}{x}e^{-x^2/t^2}\big(1+o(1)\big),\qquad t\downarrow0,
\end{equation*}
uniformly for $x\in I$. It is nonzero for every sufficiently small positive time and flat to every order at zero. The common interface pressure is $P_*+\veps e^{-x^2/t^2}$ for $t>0$ and $P_*$ at $t=0$.
\end{thm}

Setting $\varepsilon=0$ in the construction underlying Theorem 2.3 yields a reference solution with the interface $y=bt^2/(2x)$ and a constant pressure trace $P_*$. This reference solution remains an entropy wave throughout the local time interval. For $\varepsilon>0$, the constructed solutions have the same interface and satisfy initial entropy-wave material compatibility at every order, yet develop nonzero tangential slip at every sufficiently small positive time. Thus, matching the initial material jets of pressure and velocity across the interface at every order does not ensure persistence of the entropy-wave property in the smooth category. This construction has its own reference wave, distinct from the planar reference wave in Section 2.2.




\subsection{Main ideas of the proof}
The two constructions begin with a prescribed interface and a common pressure trace. The kinematic and Bernoulli relations determine Cauchy data for the phasewise velocity potentials. We then construct bulk solutions realizing these traces. In this way, the bulk initial data are obtained together with the contact history.

For Theorem~\ref{thm:main}, we prescribe the accelerating flat interface $f(t)=bt^2/2$ and the analytic pressure trace $q=P_*+\veps t^\ka a(x)$. In the accelerating coordinates, the tangential velocity potentials satisfy Hamilton--Jacobi equations on the interface. Their solutions, together with the common normal velocity, provide Cauchy data for the scalar potential equation in the normal variable. The Cauchy--Kowalevski theorem gives analytic phase fields in a collar of the interface. The different factors $R_\pm(q)^{-1}$ in the tangential pressure accelerations yield the first velocity defect at order $\ka+1$. Finally, we extract the initial fields from the collar solution, extend them smoothly, and solve two Euler Cauchy problems in physical time. Finite propagation preserves the constructed contact near the interface and gives a solution in the slab with impermeable walls.

For Theorem~\ref{thm:flat-nonlinear}, the prescribed pressure perturbation is flat at $t=0$, so a different bulk construction is needed. We use velocity potentials of the form $\phi\ppm(t,x,y)=xF\ppm(t/x,y/x)$ and the curved interface $f(t,x)=bt^2/(2x)$. This homogeneous ansatz reduces the potential equation to two independent variables. Near the prescribed interface, the reduced equation remains hyperbolic when the transverse coordinate is taken as the evolution variable, allowing smooth Cauchy data. It therefore realizes the flat pressure trace $P_*+\veps e^{-x^2/t^2}$ while preserving all initial material compatibility conditions. An equation for the difference of the boundary potentials then yields the nonzero flat slip.

The two proofs use the common formulation above but are independent. Section~\ref{sec:construction} proves Theorem~\ref{thm:main}, including the localization to fixed walls, and Section~\ref{sec:smooth} proves Theorem~\ref{thm:flat-nonlinear}.

\section{Construction with finite-order compatibility}\label{sec:construction}
We prove Theorem~\ref{thm:main} in this section.

\subsection{Construction of the interface data}\label{subsec:interface-data}
Fix the parameters in Theorem~\ref{thm:main} and prescribe
\begin{equation*}
f(t)=\tfrac12bt^2,\qquad q(t,x)=P_*+\veps t^\ka a(x).
\end{equation*}
These are the interface position and common pressure trace for a short time interval. The kinematic condition fixes the normal velocity at $f_t=bt$. To construct the tangential velocity traces, for each sign solve for $A_\pm(t,x)$:
\begin{equation}\label{eq:HJ}
\p_tA_\pm+\tfrac12|\nabla_xA_\pm|^2+H_\pm(q)=0,
\qquad A_\pm(0,x)=0.
\end{equation}
The analytic Cauchy theorem gives a solution for small $|t|$, uniformly on the compact torus. Including $\veps$ as an additional variable on which no differentiation is imposed gives joint analyticity near $\veps=0$.

This step produces the two interface potentials $A_\pm$ and the tangential velocities $v_\pm=\nabla_xA_\pm$. Together with $V\ppm=bt$, $p\ppm=q$, and $\rho\ppm=R_\pm(q)$ on the interface, these are the traces to be realized by the bulk fields.

\subsection{Construction of the bulk fields near the interface}\label{subsec:bulk-construction}
Set $\xi=y-f(t)$. With $A_\pm$ fixed, we now construct potentials $\psi\ppm(t,x,\xi)$ in a two-sided collar of $\xi=0$. We will recover the Euler fields by
\begin{equation}\label{eq:ansatz}
\ba
u\ppm&=bt\,e_y+\nabla_{x,\xi}\psi\ppm,\qquad \rho\ppm=R_\pm(p\ppm),\\
p\ppm&=H_\pm^{-1}\!\left(-\psi_t\ppm-\tfrac12|\nabla\psi\ppm|^2-b\xi\right).
\ea
\end{equation}
Here $H_\pm^{-1}$ is analytic near $H_\pm(P_*)=0$ since $H_\pm'=R_\pm^{-1}>0$. Prescribe the normal Cauchy data
\begin{equation}\label{eq:normal-data}
\psi\ppm(t,x,0)=A_\pm(t,x),\qquad
\psi_\xi\ppm(t,x,0)=0.
\end{equation}
The first datum realizes the tangential velocity $v_\pm$ and the second realizes the normal velocity $bt$. The boundary equation \eqref{eq:HJ} then makes the pressure in \eqref{eq:ansatz} equal to $q$ at $\xi=0$. All derivatives of $\psi$ are in the accelerating coordinates. The physical material derivative becomes $\p_t+\nabla\psi\cdot\nabla$ in these coordinates.

Suppress the phase label and put $v=\nabla\psi$, $h=H(p)$, and $c^2=\ga p/R(p)$. The Bernoulli relation in \eqref{eq:ansatz} gives
\begin{equation*}
(\p_t+v\cdot\nabla)h
=-\psi_{tt}-2v\cdot\nabla\psi_t-v\cdot(\nabla^2\psi)v-b\psi_\xi.
\end{equation*}
Since $h_p=\rho^{-1}$, the pressure equation is equivalent to $D_t h+c^2\dvg v=0$. Hence the potential equation is
\begin{equation}\label{eq:potential}
\psi_{tt}+2\nabla\psi\cdot\nabla\psi_t
+\nabla\psi\cdot(\nabla^2\psi)\nabla\psi
-c^2\Lap\psi+b\psi_\xi=0.
\end{equation}
The sound speed is the analytic function
\begin{equation}\label{eq:sound}
c^2=c_0^2-(\ga-1)
\left(\psi_t+\tfrac12|\nabla\psi|^2+b\xi\right),
\qquad c_0^2=\frac{\ga P_*}{r}.
\end{equation}

More explicitly, write $w=\psi_\xi$ and $z=\nabla_x\psi$. Equation \eqref{eq:potential} has the normal form
\begin{equation}\label{eq:normal-form}
\ba
(c^2-w^2)\psi_{\xi\xi}
={}&\psi_{tt}+2w\psi_{\xi t}+2z\cdot\nabla_x\psi_t
+2wz\cdot\nabla_x\psi_\xi\\
&+z\cdot(\nabla_x^2\psi)z-c^2\Lap_x\psi+bw.
\ea
\end{equation}
On $\xi=0$, $w=0$ and $c^2=\ga q/R(q)>0$. Thus the right side divided by $c^2-w^2$ is analytic in its arguments near the data. The second-order Cauchy--Kowalevski theorem applies to \eqref{eq:normal-data}--\eqref{eq:normal-form}, with $\xi$ as the Cauchy variable and $(t,x)$ as coordinates along the Cauchy surface. This theorem permits mixed derivatives with fewer than two normal differentiations and tangential derivatives of total order at most two; see \cite{Evans} for the analytic Cauchy theorem.

Local solutions around points of the torus agree on overlaps by analytic uniqueness. Compactness gives a common collar. Applying the same argument with $\veps$ among the independent analytic variables gives joint analyticity on
\begin{equation*}
|t|<T_1,\qquad |\xi|<\de_1,\qquad |\veps|<\veps_1,
\qquad x\in\T^{d-1}.
\end{equation*}
We subsequently work on a smaller closed product neighborhood. Positivity of pressure, density, and $c^2-w^2$ follows by shrinking it. This step has constructed two potentials $\psi\ppm$, each defined on both sides of $\xi=0$ with its own phase law. Formula \eqref{eq:ansatz} now defines the corresponding auxiliary Euler fields throughout the collar. To obtain the physical fields, we take the plus field on $\xi>0$ and the minus field on $\xi<0$. We next verify that these fields solve the Euler equations and form the prescribed contact discontinuity.

\subsection{Verification of the Euler equations and contact conditions}\label{subsec:contact-verification}
Taking the gradient of Bernoulli gives
\begin{equation*}
\p_t v+(v\cdot\nabla)v+\rho^{-1}\nabla p+b e_y=0,
\end{equation*}
because $v$ is a gradient. In physical coordinates $D_t u=b e_y+(\p_t+v\cdot\nabla)v$, so the inertial term cancels and the physical momentum equation in \eqref{eq:euler-conservative} is unforced. Equation \eqref{eq:potential} supplies the pressure equation. Finally,
\begin{equation*}
R'(p)=\frac{R(p)}{\ga p}
\end{equation*}
implies the mass equation. Since $p=K\rho^\ga$, entropy transport also holds. This verifies the full, rather than merely barotropic momentum, system.

At $\xi=0$, equations \eqref{eq:HJ} and \eqref{eq:ansatz} give
\begin{equation}\label{eq:traces}
p\ppm=q,\qquad V\ppm=bt=f_t,\qquad
U\ppm=v_\pm:=\nabla_x A_\pm.
\end{equation}
The interface is therefore material for both phases and satisfies \eqref{eq:contact}.

For completeness, the flux check can be written directly. The relative mass flux is $\rho(V-f_t)=0$. The relative momentum flux is $\rho u(V-f_t)+p e_y=p e_y$, whose jump vanishes. The relative energy flux is $E(V-f_t)+pV=p f_t$, whose jump also vanishes. Thus the piecewise fields solve the conservative Euler equations across the contact.

Differentiating Bernoulli normally on the interface, and using $\psi_\xi=\psi_{\xi t}=\nabla_x\psi_\xi=0$, yields
\begin{equation*}
R_\pm(q)^{-1}\p_\xi p\ppm=-b.
\end{equation*}
Consequently \eqref{eq:sign-exact} holds. Since
\begin{equation*}
R_-(q)-R_+(q)=(r_--r_+)(q/P_*)^{1/\ga},
\end{equation*}
a uniform lower bound for $q$ proves the strict sign and nonvanishing density jump.

\subsection{Verification of compatibility and vortex-sheet generation}\label{subsec:compatibility-verification}
We now verify the initial material compatibility and compute the tangential slip of the constructed fields. Tangential differentiation of \eqref{eq:HJ} gives
\begin{equation}\label{eq:burgers}
\p_t v_\pm+(v_\pm\cdot\nabla_x)v_\pm
=-\frac{\veps t^\ka}{R_\pm(q)}\nabla_x a,
\qquad v_\pm(0,x)=0.
\end{equation}
Taylor recursion in $t$ shows
\begin{equation}\label{eq:v-factor}
v_\pm(t,x,\veps)
=-\frac{\veps t^{\ka+1}}{(\ka+1)r_\pm}\nabla_x a(x)
+\veps t^{\ka+2}B_\pm(t,x,\veps),
\end{equation}
where $B_\pm$ are analytic on a smaller common neighborhood. To justify the factor $\veps$ in the remainder, observe that $A_\pm$ and $v_\pm$ vanish identically when $\veps=0$. The coefficients through degree $\ka$ vanish for every $\veps$, and the coefficient of degree $\ka+1$ is exactly the one displayed. Analytic divisibility by $\veps t^{\ka+2}$ follows. This also gives uniform bounds after any fixed number of tangential derivatives.

By \eqref{eq:trace-material}, the trace material operators are
\begin{equation*}
L_\pm=\p_t+v_\pm\cdot\nabla_x.
\end{equation*}
Since $v_\pm=O(t^{\ka+1})$, expanding a product of these operators gives
\begin{equation}\label{eq:L-jets}
(L_\pm^jG)(0,x)=(\p_t^jG)(0,x),\qquad 0\leq j\leq \ka+1,
\end{equation}
for every smooth trace $G$. Indeed, a term different from $\p_t^jG$ contains a time derivative of $v_\pm$ of order at most $j-1\leq \ka$, possibly with tangential derivatives, and vanishes at zero.

Applying \eqref{eq:L-jets} to \eqref{eq:traces} gives, on the initial interface,
\begin{equation*}
u_0\ppm=0,\qquad (D_t u\ppm)_0=b e_y,\qquad
(D_t^j u\ppm)_0=0\quad(2\leq j\leq \ka),
\end{equation*}
and
\begin{equation*}
p_0\ppm=P_*,\qquad
(D_t^jp\ppm)_0=0\quad(1\leq j<\ka),\qquad
(D_t^\ka p\ppm)_0=\veps \ka!a.
\end{equation*}
Empty index ranges are omitted when $\ka=1$. These are interface identities, not identities throughout the initial phases. At the next velocity order \eqref{eq:v-factor} gives \eqref{eq:first-defect}. Subtracting the two expansions gives \eqref{eq:slip}. If
\begin{equation*}
A_s=\frac1{\ka+1}\left(\frac1{r_+}-\frac1{r_-}\right)
\norm{\nabla_x a}_{H^s}>0,
\end{equation*}
then the reverse triangle inequality yields
\begin{equation*}
\norm{\jump U(t)}_{H^s}\geq |\veps|t^{\ka+1}(A_s-C_st).
\end{equation*}
Taking $T\leq A_s/(2C_s)$ proves \eqref{eq:lower-bound} with $c_s=A_s/2$.

\subsection{Construction of the solution in the slab}\label{subsec:localization}
The preceding steps give a contact solution in a collar of the interface. We now construct a solution in $\T^{d-1}\times(-L,L)$ with impermeable walls. Since $f(0)=0$, evaluation of the auxiliary fields at $t=0$ gives initial data on a two-sided collar of $y=0$:
\begin{equation*}
\cU_{0,\veps}\ppm(x,y)=\cU_\veps\ppm(0,x,y),\qquad \cU=(\rho,u,p).
\end{equation*}
These are the bulk initial data produced by the normal Cauchy construction. We will extend them away from the interface and use them as Cauchy data for a new evolution in physical time.

First, we record their dependence on $\veps$. At $\veps=0$, equation \eqref{eq:HJ} gives $A_\pm=0$. Analytic uniqueness for the normal problem gives $\psi\ppm=0$, and \eqref{eq:ansatz} recovers \eqref{eq:background} with physical velocity $bt\,e_y$. On a compact smaller collar, joint analyticity gives
\begin{equation*}
\sup_{|\veps|\leq\veps_0}
\norm{\p_\veps\cU_{0,\veps}\ppm}_{C^m}<\infty.
\end{equation*}
Integration in $\veps$ proves \eqref{eq:smallness}. The same argument supplies uniform bounds of every fixed order needed below.

Next, choose a smooth cutoff $\chi(y)$ supported strictly inside the initial collar and $(-L,L)$, and equal to one on a smaller collar. Choose constants $\rho_{*,\pm}>0$, and extend each auxiliary initial field to $\T^{d-1}\times\R$ by
\begin{equation*}
\rho_{0,\veps}^{\mathrm{ext},\pm}
=\rho_{*,\pm}+\chi(y)\big(\rho_{0,\veps}\ppm-\rho_{*,\pm}\big),\qquad
u_{0,\veps}^{\mathrm{ext},\pm}=\chi(y)u_{0,\veps}\ppm,
\end{equation*}
with $0\leq\chi\leq1$, and set $p_{0,\veps}^{\mathrm{ext},\pm}=K_\pm(\rho_{0,\veps}^{\mathrm{ext},\pm})^\ga$. The products are defined using the collar fields only on a neighborhood of the support of $\chi$. Positivity is preserved. The same cutoff and constants are used for all $\veps$, including zero, so the full initial Sobolev difference is $O(|\veps|)$.

Finally, solve the two smooth Euler Cauchy problems independently on the whole cylinder with these extended initial data, and denote their solutions by $\cU_\veps^{\mathrm{ext},\pm}$. The local existence theorem for quasilinear symmetric hyperbolic systems \cite{Kato1975,VH1972} gives a common short lifespan from the uniform positive lower bounds and a fixed uniform high Sobolev bound on the deviations from the phasewise constant equilibria. Smooth initial data remain smooth on that lifespan.

Shorten time so that the backward domains of dependence of a smaller moving collar meet $t=0$ entirely where $\chi=1$, while disturbances from the cutoff support cannot reach the walls. Finite propagation and local uniqueness identify $\cU_\veps^{\mathrm{ext},\pm}$ with the corresponding collar fields there. Near $y=\pm L$, each extended solution remains its constant equilibrium. Restricting $\cU_\veps^{\mathrm{ext},+}$ to $f(t)<y<L$ and $\cU_\veps^{\mathrm{ext},-}$ to $-L<y<f(t)$ gives the asserted contact in the full slab, with $V=0$ at the fixed walls. Since it agrees with the collar construction near the interface, all contact conditions, compatibility identities, and slip estimates proved above persist. This completes the proof of Theorem~\ref{thm:main}.

\section{Smooth infinite-order compatibility and sheet generation}\label{sec:smooth}
We prove Theorem~\ref{thm:flat-nonlinear} in this section.

\subsection{A homogeneous reduction of the potential equation}
Work in two spatial dimensions on a tangential interval whose closure lies in $(0,\infty)$. The homogeneous ansatz reduces the transverse Cauchy problem to an equation that is hyperbolic near the reference state, allowing us to realize the prescribed smooth, nonanalytic interface data. We seek the phasewise velocity potentials in the form
\begin{equation}\label{eq:homogeneous}
s=\frac{t}{x},\qquad z=\frac{y}{x},\qquad
\phi\ppm(t,x,y)=xF\ppm(s,z).
\end{equation}
Here $F^\pm$ are the unknown reduced potentials, and $s$ denotes a similarity coordinate, not the entropy. For one phase, suppress the sign. Once $F$ is constructed, the velocity, pressure, and density are recovered by
\begin{equation}\label{eq:homogeneous-velocity}
\ba
V&=\phi_y=F_z,\qquad U=\phi_x=F-sF_s-zF_z,\\
p&=H^{-1}\!\left(-F_s-\tfrac12(U^2+V^2)\right),\qquad \rho=R(p),
\ea
\end{equation}
where $R,H$ are defined in \eqref{eq:R-H}, and $H^{-1}$ is the local inverse near $H(P_*)=0$. Thus $u=(U,V)=\nabla_{x,y}\phi$. The sound speed satisfies
\begin{equation*}
c^2=\frac{\ga p}{\rho}
=\frac{\ga P_*}{r}-(\ga-1)\big(F_s+\tfrac12(U^2+V^2)\big).
\end{equation*}
Bernoulli's identity gives momentum. Mass conservation is equivalent to the potential equation
\begin{equation}\label{eq:physical-potential}
\phi_{tt}+2\nabla\phi\cdot\nabla\phi_t
+\nabla\phi\cdot D^2\phi\nabla\phi-c^2\Lap\phi=0.
\end{equation}
The pressure equation then follows from $p=K\rho^\ga$. Substitution of \eqref{eq:homogeneous} gives the exact reduced equation
\begin{equation}\label{eq:reduced}
\mathsf A F_{ss}+2\mathsf B F_{sz}+\mathsf C F_{zz}=0,
\end{equation}
where
\begin{equation}\label{eq:reduced-coefficients}
\ba
\mathsf A&=(1-sU)^2-c^2s^2,\\
\mathsf B&=(1-sU)(V-zU)-c^2sz,\\
\mathsf C&=(V-zU)^2-c^2(1+z^2).
\ea
\end{equation}
For verification, the independent second derivatives, multiplied by $x$, are
\begin{equation*}
\ba
x\phi_{tt}&=F_{ss},\quad x\phi_{tx}=-sF_{ss}-zF_{sz},\quad
x\phi_{ty}=F_{sz},\\
x\phi_{xx}&=s^2F_{ss}+2szF_{sz}+z^2F_{zz},\quad
x\phi_{xy}=-sF_{sz}-zF_{zz},\\
x\phi_{yy}&=F_{zz}.
\ea
\end{equation*}
Thus it remains to construct $F^\pm$ satisfying \eqref{eq:reduced} with Cauchy data that enforce the contact conditions.

\subsection{Construction of the interface data}
We first prescribe the interface and its pressure, and construct the potential traces that will serve as Cauchy data for \eqref{eq:reduced}. Fix $b>0$, $P_*>0$, $r_->r_+>0$, and $\veps>0$. Define
\begin{equation}\label{eq:flat-forcing}
\chi(s)=\begin{cases}e^{-1/s^2},&s>0,\\0,&s\leq0,\end{cases}
\qquad Q(s)=P_*+\veps\chi(s),\qquad g(s)=\tfrac12bs^2.
\end{equation}
We prescribe $Q$ as the common pressure on the reduced curve $z=g(s)$, which corresponds to the physical interface
\begin{equation}\label{eq:flat-interface}
y=f(t,x)=x g(t/x)=\frac{bt^2}{2x}.
\end{equation}
Let $a_\pm(s)=F\ppm(s,g(s))$ denote the unknown potential traces on this curve. Since $s=t/x$, differentiating $f(t,x)=xg(s)$ gives
\begin{equation*}
f_t=g'(s),\qquad f_x=g(s)-sg'(s)=:\vartheta(s)=-\tfrac12bs^2.
\end{equation*}
By \eqref{eq:homogeneous}, the required potential traces in physical variables are
\begin{equation*}
\phi\ppm(t,x,f(t,x))=x a_\pm(t/x).
\end{equation*}
Write $(U_\pm,V_\pm)$ for the boundary velocities. Differentiating this trace identity with respect to $x$ at fixed $t$, and using $u\ppm=\nabla_{x,y}\phi\ppm$, yields
\begin{equation*}
U_\pm+\vartheta V_\pm
=\p_x\big(x a_\pm(t/x)\big)=a_\pm-sa_\pm'.
\end{equation*}
On the other hand, the kinematic condition $f_t+U_\pm f_x=V_\pm$ becomes
\begin{equation*}
V_\pm-\vartheta U_\pm=g'.
\end{equation*}
Solving these two linear equations for $U_\pm$ and $V_\pm$ gives
\begin{equation}\label{eq:flat-boundary-velocities}
U_\pm=\frac{a_\pm-sa_\pm'-\vartheta g'}{1+\vartheta^2},\qquad
V_\pm=\frac{g'+\vartheta(a_\pm-sa_\pm')}{1+\vartheta^2}.
\end{equation}
Differentiating $a_\pm(s)=F\ppm(s,g(s))$ with respect to $s$ also gives $a_\pm'=F_s\ppm+g'F_z\ppm$ on the curve. Hence $\phi_t\ppm=F_s\ppm=a_\pm'-g'V_\pm$ on the interface, and Bernoulli's relation with $p\ppm=Q$ gives the scalar implicit ODE
\begin{equation}\label{eq:flat-ode}
a_\pm'-g'V_\pm+\tfrac12(U_\pm^2+V_\pm^2)
+H_\pm(Q)=0,\qquad a_\pm(0)=0,
\end{equation}
where $U_\pm,V_\pm$ are the expressions in \eqref{eq:flat-boundary-velocities}. The derivative of the left side with respect to $a_\pm'$ is $1-sU_\pm$, which equals one at $(s,a_\pm,a_\pm')=(0,0,0)$. The implicit-function theorem therefore solves \eqref{eq:flat-ode} for $a_\pm'$ as a smooth function of $(s,a_\pm)$. The smooth ODE theorem supplies unique solutions on an interval about zero.

The functions $a_\pm$ and the boundary velocities \eqref{eq:flat-boundary-velocities} are now fixed. They supply the Cauchy data
\begin{equation}\label{eq:flat-cauchy}
F\ppm(s,g(s))=a_\pm(s),\qquad
F_z\ppm(s,g(s))=V_\pm(s).
\end{equation}

\subsection{Construction of the bulk fields near the interface}
We now solve \eqref{eq:reduced} for $F^\pm$ on both sides of $z=g(s)$ with the data \eqref{eq:flat-cauchy}. These data give $F=F_s=F_z=0$ at the origin, and hence
\begin{equation*}
\mathsf A=1,\qquad \mathsf B=0,\qquad \mathsf C=-\frac{\ga P_*}{r_\pm}<0.
\end{equation*}
Flatten the curve by $\eta=z-g(s)$ and set $\wt F(s,\eta)=F(s,\eta+g(s))$. The coefficients of $\wt F_{ss}$, $2\wt F_{s\eta}$, and $\wt F_{\eta\eta}$ become respectively
\begin{equation*}
\mathsf A,\qquad
\wt{\mathsf B}=\mathsf B-g'\mathsf A,\qquad
\wt{\mathsf C}=\mathsf C-2g'\mathsf B+(g')^2\mathsf A.
\end{equation*}
In these coordinates,
\begin{equation*}
U=\wt F-s\wt F_s-(\eta+\vartheta)\wt F_\eta,\qquad
V=\wt F_\eta,
\end{equation*}
and the exact flattened equation, including its curvature term, is
\begin{equation}\label{eq:flattened-reduced}
\mathsf A\wt F_{ss}
+2\wt{\mathsf B}\wt F_{s\eta}
+\wt{\mathsf C}\wt F_{\eta\eta}
-\mathsf A g''\wt F_\eta=0.
\end{equation}
All coefficients are smooth functions of $(s,\eta,\wt F,\wt F_s,\wt F_\eta)$.
Near the origin $\mathsf A>0$ and $\wt{\mathsf C}<0$. The first derivatives satisfy
\begin{equation}\label{eq:reduced-first-order}
\p_\eta
\begin{pmatrix}\wt F_\eta\\\wt F_s\end{pmatrix}
+
\begin{pmatrix}
2\wt{\mathsf B}/\wt{\mathsf C}
&\mathsf A/\wt{\mathsf C}\\
-1&0
\end{pmatrix}
\p_s\begin{pmatrix}\wt F_\eta\\\wt F_s\end{pmatrix}
=\begin{pmatrix}(\mathsf A/\wt{\mathsf C})g''\wt F_\eta\\0\end{pmatrix}.
\end{equation}
For the first-order Cauchy problem, regard the two displayed derivative components as independent unknowns and adjoin the equation that the $\eta$-derivative of $\wt F$ equals the first component. The data on $\eta=0$ are
\begin{equation*}
(\wt F,\wt F_\eta,\wt F_s)=(a_\pm,V_\pm,a_\pm').
\end{equation*}
The displayed coefficient matrix has real, distinct eigenvalues since
$\wt{\mathsf B}^{2}-\mathsf A\wt{\mathsf C}>0$.
The positive matrix $\operatorname{diag}(1,-\mathsf A/\wt{\mathsf C})$ symmetrizes it; adjoining the scalar $\wt F$ equation preserves symmetric hyperbolicity. The second row and the equation for $\wt F$ preserve the identity of the second component with $\p_s\wt F$, which holds at $\eta=0$. Thus the first-order solution reconstructs a solution of \eqref{eq:flattened-reduced}. The smooth local existence and uniqueness theorem for quasilinear symmetric hyperbolic systems \cite{Kato1975,VH1972}, with $\eta$ as the evolution variable, applies for both signs of $\eta$. It yields two $C^\infty$ potentials $F^\pm$, each defined on a neighborhood of the reduced origin crossing $z=g(s)$.

For $x$ in any fixed compact interval $I\subset(0,\infty)$, choose $T,\de>0$ so that $|t|\leq T$, $|y-f(t,x)|<\de$, $x\in I$ map into this reduced neighborhood. Define $\phi^\pm=xF^\pm(t/x,y/x)$ and recover $u^\pm,p^\pm,\rho^\pm$ from \eqref{eq:homogeneous-velocity}, with $S^\pm=\log(K_\pm/A)$. By the reduction above, these are smooth auxiliary Euler solutions on both sides of the physical interface. Pressure, density, and sound speed are positive after shrinking the collar.

Restrict the plus and minus solutions to $y>f$ and $y<f$, respectively. Their traces are the interface data constructed above: \eqref{eq:flat-boundary-velocities} gives the kinematic condition, and \eqref{eq:flat-ode} gives the common pressure $Q(t/x)$. Their density traces are $R_\pm(Q(t/x))$, hence have a nonzero jump. This constructs the contact solution. Its initial data are obtained by evaluating the bulk fields at $t=0$, with $f_0=0$.

\subsection{Verification of compatibility and the Taylor jump}
We verify that the initial data just constructed satisfy the entropy-wave compatibility conditions of every order and the positive Taylor-jump condition. For $s\leq0$ the forcing has $Q=P_*$ and $H_\pm(Q)=0$. Thus the two boundary ODEs coincide there and have the same initial value. Consequently $a_+=a_-$ for $s\leq0$, and their difference is flat at zero. Formula \eqref{eq:flat-boundary-velocities} shows that both velocity differences are flat as well. The pressure trace $Q-P_*$ is flat. On every compact $I\subset(0,\infty)$ these properties hold jointly in $(t,x)$ after composition with $s=t/x$.

By \eqref{eq:trace-material}, on the interface the phasewise material derivatives are
\begin{equation*}
\p_t+U_\pm(t/x)\p_x.
\end{equation*}
The traces and the coefficients of these operators have identical time jets at zero. Repeated differentiation therefore gives
$\jump{(D_t^ju)_0}=\jump{(D_t^jp)_0}=0$ for every $j\geq0$.
These are actual jets of the constructed Euler solution, hence also the jets obtained recursively from its initial data.

The ODE gives $a_\pm'(0)=a_\pm''(0)=0$. It follows that the initial interface velocity is zero and
\begin{equation*}
\left.\frac{\mathrm d}{\mathrm dt}(U_\pm,V_\pm)(t/x)\right|_{t=0}
=\left(0,\frac{b}{x}\right).
\end{equation*}
Initially $f_t=f_x=0$, $\rho\ppm=r_\pm$, and momentum gives
\begin{equation}\label{eq:flat-taylor}
\p_yp_0\ppm(x,0)=-\frac{br_\pm}{x},\qquad
\cT(0,x)=\frac{b(r_--r_+)}{x}>0.
\end{equation}
Compactness of $I$ and continuity preserve a uniform positive Taylor jump for a sufficiently short time.

\subsection{Verification of vortex-sheet generation}
It remains to show that the constructed contact solution has nonzero tangential slip for every sufficiently small positive time. Set
\begin{equation*}
Z=a_+-a_-,\qquad \Delta H=H_+(Q)-H_-(Q),\qquad
m=\tfrac12(U_++U_-).
\end{equation*}
Subtracting the two equations in \eqref{eq:flat-ode} and using $V_\pm-g'=\vartheta U_\pm$ yields
\begin{equation}\label{eq:flat-difference}
(1-sm)Z'+mZ=-\Delta H,
\qquad Z(0)=0.
\end{equation}
For small positive $s$, $1-sm>0$ and
\begin{equation*}
\Delta H=\veps\left(\frac1{r_+}-\frac1{r_-}\right)
\chi(s)\big(1+O(\veps\chi(s))\big)>0.
\end{equation*}
The integrating-factor formula in \eqref{eq:flat-difference} shows that $Z<0$ and
\begin{equation*}
|Z(s)|\leq C\veps\int_0^s\chi(\sig)\,\rd\sig
\leq \tfrac12 C\veps s^3\chi(s),
\end{equation*}
where $C$ may depend on the fixed densities $r_\pm$. The last inequality follows by integrating
$\chi'=2s^{-3}\chi$; explicitly,
\begin{equation*}
\int_0^s\chi(\sig)\,\rd\sig
=\tfrac12s^3\chi(s)-\tfrac32\int_0^s\sig^2\chi(\sig)\,\rd\sig.
\end{equation*}
With unit tangent $\tau=(1,\vartheta)/\sqrt{1+\vartheta^2}$, the exact sheet strength is
\begin{equation}\label{eq:flat-slip-exact}
\jump{u_\tau}
=\frac{Z-sZ'}{\sqrt{1+\vartheta^2}}
=\frac{Z+s\Delta H}{(1-sm)\sqrt{1+\vartheta^2}}.
\end{equation}
Since $Z/(s\Delta H)=O(s^2)$ and $m,\vartheta\to0$, this proves
\begin{equation}\label{eq:flat-slip-asymptotic}
\jump{u_\tau}\big(t,x,f(t,x)\big)
=\veps\left(\frac1{r_+}-\frac1{r_-}\right)
\frac{t}{x}e^{-x^2/t^2}\big(1+o(1)\big),\qquad t\downarrow0,
\end{equation}
uniformly on $I$. In particular the jump is positive at every point of $I$ for every sufficiently small positive time, despite its zero jet of every order at zero. This completes the proof of Theorem~\ref{thm:flat-nonlinear}.

\textbf{Acknowledgments.}
Yuri Trakhinin is supported by the state assignment of the Sobolev Institute of Mathematics (project No. FWNF-2026-0028).
Wenbin Zhao is partially supported by the National Natural Science
Foundation of China (No.\ 12401286).

\textbf{Data availability statement.}
There are no data associated with this article.

\textbf{Conflict of interest.}
The authors declare that they have no conflict of interest.

\textbf{AI assistance statement.} The authors used AI tools to assist with exploring and
manuscript editing; all mathematical validation, final proof decisions, and final wording
remain the sole responsibility of the human authors.

\end{document}